# The Reconstruction of Graphs

Dhananjay P. Mehendale
Sir Parashurambhau College, Tilak Road, Pune-411030,
India.

## Abstract

In this paper we discuss reconstruction problems for graphs. We develop some new ideas like isomorphic extension of isomorphic graphs, partitioning of vertex sets into sets of equivalent points, subdeck property, etc. and develop an approach to deal with reconstruction problem. We then discuss complete sets of invariants for graphs and reconstruction conjecture. We then begin with development of few equivalent formulations of reconstruction conjecture. In the last section we briefly elaborate the formulation due to Harary its exact demand and finally proceed to give a different proof of reconstruction conjecture using reconstructibility of graph from its spanning trees and reconstructibility of tree from its pendant point deleted deck of subtrees. This last proof can be used to develop a systematic procedure to reconstruct unique graph from its deck.

**1. Introduction:** This well known conjecture, called the reconstruction conjecture, of Kelly [1] and Ulam [2] can be stated as follows: Let $G$ and $H$ be two simple graphs with p points, p>2, and let there exist a (1-1), onto map,

$\theta: V(G) \to V(H)$,
$\quad v_i \in V(G) \to u_i \in V(H)$

such that $G_i \equiv G - v_i \cong H_i \equiv H - u_i, \forall i, i = 1, 2, \cdots, p$, then
$G \cong H$. (The symbol $\cong$ stands for isomorphism.)

The existence of such a map is called hypomorphism. The reconstruction conjecture states that hypomorphism implies isomorphism.

Let $v_i \in V(G)$. Then the induced subgraph $G_i \equiv G - v_i$ is called the vertex deleted subgraph (VDS) of $G$. The collection of all VDS, $\{G_i \equiv G - v_i, i = 1, 2, \cdots, p\}$, is called the deck of (vertex deleted subgraphs of) $G$, $Deck(G)$. The graph $G$ is reconstructible if every graph with the same deck as $G$ is isomorphic to $G$. Thus, the conjecture can be stated as: Every graph with more than two vertices is uniquely (up to isomorphism) reconstructible from its deck.

In order to reconstruct $G$ from its deck, $Deck(G)$, it suffices to identify which must be the neighbors (adjacent vertices) of the missing vertex $v_i$ in some graph $G_i$ in the deck [3]. We proceed to discuss about how to accomplish this task.

**2. Preliminaries:** Let $G$ and $H$ be two simple (p, q) graphs, i.e. having p points (vertices) and q lines (edges), and let $G$ be isomorphic to H. Then there exists a (1-1), onto, adjacency preserving map

$$\theta: V(G) \to V(H),$$
$$v_i \in V(G) \to u_i \in V(H)$$

such that if a pair of vertices $(v_k, v_l) \in V(G)$, the vertex set for $G$, are adjacent (nonadjacent) then the corresponding pair of vertices, under the isomorphism map $\theta$, the pair of vertices $(u_k, u_l) \in V(H)$, the vertex set for $H$, will be adjacent (nonadjacent).

Now, if we add a new vertex $v$ to the vertex set of $G$, $V(G)$, and a new vertex $u$ to the vertex set of $H$, $V(H)$, and also add an edge $(v, v_i)$ to the edge set of $G$, $E(G)$, and the corresponding edge $(u, u_i)$ to the edge set of $H$, $E(H)$, and denote the extended $G$ and $H$ by $G^e$ and $H^e$ respectively, then clearly $G^e$ and $H^e$ will be isomorphic with the isomorphism map

$$\theta^e : V(G^e) \to V(H^e) \text{ where}$$
$$\theta^e(V(G)) = \theta(V(G)), \text{ and}$$
$$\theta^e(v) = u.$$

Few definitions are now in order:

**Definition 2.1:** Let $G$ and $H$ be two isomorphic graphs. The process described above of extending graphs $G$ and $H$ to their respective supergraphs $G^e$ and $H^e$ is called the **isomorphic extension of isomorphic graphs.**

**Definition 2.2:** Let $G$ be an unlabelled (p, q) graph and let $G^e$ be a supergraph of $G$ obtained by taking a (new) vertex outside of the vertex set $V(G)$ and joining it to some (unspecified) vertex of $G$ by an edge not in the edge set $E(G)$ is called the **extension** of $G$ to $G^e$.

**Definition 2.3.1:** A subset $V_j$ of points $\{u_1^j, u_2^j, \cdots, u_r^j\}$ in a tree T is called a **set of equivalent points** or simply an **equivalent set** if all the trees $T + vu_s^j$, $1 \leq s \leq r$, obtained from T by adding an edge $vu_s^j$, obtained by joining vertex $u_s^j$ in set $V_j$ to a new vertex $v$, are all isomorphic to each other.

**Definition 2.3.2:** A subset of points (vertices) of V(G) is called a **set of equivalent points** or an **equivalent set** if the graphs obtained by extension of G at any point among these points, achieved by joining any point among these points, one at a time, to a (new) point outside, not in V(G), leads to graphs which are all mutually isomorphic.

**Definition 2.4:** A collection of subsets $\{V_1, V_2, \cdots, V_m\}$ of V(G), the vertex set of G is called **partitioning of V(G) into equivalent sets** if all subsets $V_i$, i = 1, 2,…, m are equivalent sets,

$$V_i \cap V_j = \phi, \forall i \neq j, \text{ and } V(G) = \bigcup_{i=1}^{m} V_i$$, where $\phi$ is a null set.

**Theorem 2.1:** Let $G$ be a (p, q) graph and let $\{V_1, V_2, \cdots, V_k\}$ be the partitioning of the vertex set of $G$ into equivalent sets, and let $\{|V_1|, |V_2|, \cdots, |V_k|\}$ be their cardinalities. The number of nonisomorphic graphs that can be achieved from a graph G by adding some r edges, $r \leq p$, by extension, i.e. by joining some r points of G to an externally taken new point not in the vertex set of $G$, is equal to the number of (restricted)
k-compositions of r, like: $r = \alpha_1 + \alpha_2 + \cdots + \alpha_j + \cdots + \alpha_k$, where $0 \leq \alpha_j \leq |V_j|$.

**Proof:** Any two graphs with identical composition (though distinct ones due to extension at distinct points from the same equivalent set) are clearly isomorphic (refer to definition 2.3.2). Any two graphs with distinct composition will be clearly nonisomorphic, since they are arrived at by choosing different number of points for joining to the new point from the equivalent sets, differing at least at one place, say $\alpha_j$.

□

**Definition 2.5:** Let $G$ and $H$ be two graphs on p points. We say that $Deck(G)$ **is a subdeck of** $Deck(H)$, or, $Deck(G) \subseteq Deck(H)$ when there exists a bijection $\psi: V(G) \to V(H)$ such that
$$v_i \in V(G) \to u_i \in V(H), \text{ where } V(G), V(H) \text{ are}$$
the vertex sets of $G$ and $H$ respectively, and $G_i \equiv G - v_i$ is a subgraph of $H_i \equiv H - u_i$ for all i = 1, 2, …, p. Whenever $G_i \equiv G - v_i$ is a proper subgraph of $H_i \equiv H - u_i$ for at least one $i = 1, 2, \ldots, p$, we say that $Deck(G)$ **is a proper subdeck of** $Deck(H)$, i.e. $Deck(G) \subset Deck(H)$.

Thus, $Deck(G)$ is **not a subdeck of** $Deck(H)$ if there exists no map $\psi$ (see definition 2.5), such that the inclusion relation, namely, $Deck(G) \subseteq Deck(H)$ holds good, i.e. this relation is violated for at least one subgraph of $G$ and the corresponding subgraph of $H$.

**Definition 2.6:** The subisomorphism map for two **p point** graphs $G$ and $H$
is a bijection
$$\psi: V(G) \to V(H) \text{ where}$$
$$v_i \in V(G) \to u_i \in V(H)$$
such that if a pair of vertices $(v_k, v_l) \in V(G)$ are adjacent then the corresponding pair of vertices, under the subisomorphism map, $(u_k, u_l) \in V(H)$, will be adjacent and when $(v_k, v_l) \in V(G)$ are nonadjacent then $(u_k, u_l) \in V(H)$ **may or may not be** adjacent.

**Remark 2.1:** Under the existence of subisomorphism map $\psi: V(G) \to V(H)$ $G$ will be a spanning subgraph of $H$.

**Definition 2.7:** Checking whether $Deck(G) \subseteq Deck(H)$ holds or doesn't hold for some graphs $G$ and $H$ will be called **checking of subdeck property for graphs $G$ and $H$**

**Lemma 2.1:** Let $G_1$ and $G_2$ be two graphs on p points. If $G_1$ is subgraph of $G_2$ then $Deck(G_1) \subseteq Deck(G_2)$. Moreover, if $G_1$ is a proper subgraph of $G_2$ then $Deck(G_1) \subset Deck(G_2)$.

**Proof:** Since, $G_1$ is a subgraph of $G_2$ there will exist a subisomorphism map under which every adjacent vertex pair of $G_1$ will be adjacent in $G_2$ and $Deck(G_1) \subseteq Deck(G_2)$. When $G_1$ will be a proper subgraph of $G_2$ then there will exist a nonadjacent vertex pair in $G_1$ which will be adjacent in $G_2$ so when we consider a point deleted subgraph of $G_1$ and the corresponding point deleted subgraph of $G_2$ arrived at by deleting a vertex other than the one contained in the vertex pair that forms a nonadjacent vertex pair in $G_1$ and whose image is an adjacent vertex pair in $G_2$ and $Deck(G_1) \subset Deck(G_2)$ becomes clear.

□

**Lemma 2.2:** Let $G_i \equiv G - v_i$. We have $G^1 = G_i + (v, x_{j_1})$, where $(v, x_{j_1})$ is the edge obtained by joining the new point $v$ not in $V(G_i)$ to $x_{j_1} \in V(G_i)$. If $Deck(G^1) \subseteq Deck(G)$, then $G^1 \subseteq G$ (i.e. $G^1$ is subgraph of $G$).

**Proof:** Given $Deck(G^1) \subseteq Deck(G)$. So, for $x_{j_1}, x_{j_2}, \cdots, x_{j_{(p-1)}} \in G^1$ there will exist $v_{j_1}, v_{j_2}, \cdots, v_{j_{(p-1)}} \in G$ such that $G^1 - x_{j_r} \subseteq G - v_{j_r}$ and $G^1 - v \subseteq G - v_i$. We define the bijection

$$\psi : V(G^1) \to V(G) \text{ where}$$

$x_{j_r} \to v_{j_r}$ and $v \to v_i$. Now, if a pair of vertices $x_{j_p}, x_{j_q}$ is adjacent $G^1$ then corresponding pair $v_{j_p}, v_{j_q}$ will be adjacent $G$, since $G^1 - v \subseteq G - v_i$ (In fact, $G^1 - v = G_i = G - v_i$). Also if pair of vertices $v, x_{j_1}$ will be adjacent in $G^1$ then corresponding pair $v_i, v_{j_1}$ will be adjacent in $G$, while pairs $v, x_{j_r}$, $r \neq 1$, are nonadjacent in $G^1$ and the corresponding pairs $v_i, v_{j_r}$ may or may not be adjacent in $G$ since $G^1 - x_{j_r} \subseteq G - v_{j_r}$. So the map is a subisomorphism map.

□

**Remark 2.1:** The lemma will work well also in the condition when we have $G^r$ instead of $G^1$ obtained by joining the new vertex $v$ not in $V(G_i)$ to some r points of $G_i \equiv G - v_i$, $1 \leq r \leq d$ and d = degree of $v_i$, such that $Deck(G^r) \subseteq Deck(G)$. This will lead us to the same conclusion that $G^r \subseteq G$.

**Remark 2.2:** Let $G_{ij} = G_i - x_j$, $x_j \in V(G_i)$, and let $Deck(G^r) \subseteq Deck(G)$, where $G^r$ is obtained from $G_i$ by r extensions achieved by joining some r points of $G_i$ to the new point $v$ not in $V(G_i)$. There exist subgraphs $G_l = G - v_l$ in the $Deck(G)$ such that $G_{ij}$ will be subgraph of $G_l$ when $v_l$ in $G$ will correspond to (or stand for or represents the appearance of) $x_j$ in $G_i$. There will be a point in $G_l$ and not contained in $G_{ij}$, say $z_m$, of degree d or (d – 1) depending upon whether $z_m$ is adjacent or not adjacent to $v_l$. This point actually represents the appearance of $v_i$ in $G_l$, in such a way that $G^r - x_j$ will be subgraph of $G_l$. Since this will happen for all $x_j$, i.e. for every $x_j$ in $G_i$ there will correspond a $v_l$ defining a $G_l$ such that $G^r - x_j$ will be subgraph of $G_l$. So, we will thus have a subisomorphism map

$\psi : V(G^r) \to V(G)$ where

$x_j \to v_l$ and $v \to v_i$.

**3. Kelly-Ulam Conjecture and a Reconstruction Algorithm :** We now proceed to show that hypomorphic graphs are also isomorphic, i.e. the well-known Kelly-Ulam conjecture for graphs. The proof contains an algorithm for actual reconstruction.

**Theorem 3.1(Kelly-Ulam Conjecture):** Hypomorphic graphs are isomorphic.

**Proof 1:** Let $G$ and $H$ be two simple hypomorphic graphs on p points, p>2, i.e. let there exist a bijection,
$\quad \theta : V(G) \to V(H)$,
$\quad\quad v_i \in V(G) \to u_i \in V(H)$, where
$G_i \equiv G - v_i \cong H_i \equiv H - u_i, \forall i, i = 1, 2, \cdots, p$.
Let $\omega_i : V(G_i) \to V(H_i)$, be the corresponding isomorphism maps for i = 1, 2, …, p.
   (1) Consider some subgraph $G_i \equiv G - v_i$ of G and the corresponding isomorphic subgraph $H_i \equiv H - u_i$ of H.
   (2) Extend isomorphically the graphs $G_i$ to $G^1$ and $H_i$ to $H^1$ by adding an edge by joining some point $x_{j_1} \in V(G_i)$, the vertex set of $G_i$, to a new point $v$ not belonging to $V(G_i)$, leading to the formation of the graph $G^1$ and by adding an edge joining the corresponding point $y_{j_1} = \omega_i(x_{j_1}) \in V(H_i)$, the vertex set of $H_i$, to a new point $u$, not belonging to $V(H_i)$, leading to the formation of the graph $H^1$. (Note that there are in all (p-1) possibilities for adding edges, one in $G_i$ and the corresponding one under isomorphism map in $H_i$, since there are (p-1) vertices in both $G_i$ and $H_i$.)

(3) Check whether $Deck(G^1) \subseteq Deck(G)$ holds. If yes, then as per the lemma 2.2 $G^1$ is a proper extension of $G_i$ towards $G$, or $G^1$ is a subgraph of $G$, (and automatically $H^1$ will be a proper extension of $H_i$ towards $H$, or $H^1$ will be subgraph of $H$).

(4) Now, remove the added edge and add a new edge by joining some new point $x_{j_2} \in V(G_i)$, the vertex set of $G_i$, to a new point $v$, not belonging to $V(G_i)$, leading to the formation of the graph $G^1$ and by adding an edge joining the corresponding point $y_{j_2} = \omega_i(x_{j_2}) \in V(H_i)$, the vertex set of $H_i$, to a new point $u$, not belonging to $V(H_i)$, leading to the formation of the graph $H^1$. Go to step (3).

(5) Prepare the list of say, **proper** vertices i.e. **where the extension is proper**, namely $\{x_{j_s} \in V(G_i), s = 1, 2, \ldots, r\}$.

(6) From the pairs of proper vertices so obtained, $\binom{r}{2}$ in number, determine the list of **proper pairs** i.e. **for which the extension is proper**, i.e. when the extension is carried out by joining the pair of proper vertices (in $G_i$) and their images under $\omega_i$ (in $H_i$) to the external points, respectively to $v$ and $u$, leading to the formation of the graph $G^2$ and $H^2$, such that $Deck(G^2) \subseteq Deck(G)$, implying (remark 2.1) $G^2 \subseteq G$ (and $H^2 \subseteq H$).

(7) Form the triplets of proper vertices, $\binom{r}{3}$ in number, and their images in $H$ such that **any two of them also form a proper pair**, determine the **proper triplets** for which $Deck(G^3) \subseteq Deck(G)$ (and $Deck(H^3) \subseteq Deck(H)$) implying $G^3 \subseteq G$ (and $H^3 \subseteq H$).

(8) Continue extending isomorphically, and each time taking the proper action implied by the correctness (or incorrectness) of adding the edges as per the above given steps, till we reach the graphs $G^d$ and $H^d$, where d is the degrees of vertices $v_i, u_i$ (the number of vertices, number of edges, and degree sequences of hypomorphic graphs are identical). At this stage we will not only have $Deck(G^d) \subseteq Deck(G)$ and $Deck(H^d) \subseteq Deck(H)$, moreover $Deck(G^d) = Deck(G)$ implying $G^d \cong G$ and $Deck(H^d) = Deck(H)$ implying $H^d \cong H$. And by the isomorphic extension of $G_i$ to $G^d$ and $H_i$ to $H^d$, i.e. by $G^d \cong H^d$, we have $G \cong H$. □

**Proof 2:** Since the number of vertices, number of edges, degree sequences etc. of hypomorphic graphs are identical. We can take the degree of vertex $v_i, d(v_i)$ and the degree of vertex $u_i, d(u_i)$ as identical and equal to d say.

(1) Consider some subgraph $G_i \equiv G - v_i$ of G and the corresponding isomorphic subgraph $H_i \equiv H - u_i$ of H.

(2) Extend isomorphically the graphs $G_i$ to $G^d$ and $H_i$ to $H^d$ by adding some d edges by joining some d points $x_{j_1}, x_{j_2}, \cdots, x_{j_d} \in V(G_i)$, the vertex set of $G_i$, to a new point $v$, not belonging to $V(G_i)$, leading to the formation of the graph $G^d$ and by adding some d edge joining the corresponding d points $y_{j_1}, y_{j_2}, \cdots, y_{j_d} \in V(H_i)$, the vertex set of $H_i$, such that $y_{j_r} = \omega_i(x_{j_r}), r = 1, 2, \cdots, d$, to a new point $u$, not belonging to $V(H_i)$, leading to the formation of the graph $H^d$. Note that the total number of possibilities for adding d edges at the same time by joining vertex $v$ and some d points of $G_i$ and for adding d edges by joining vertex $u$ and the corresponding d points of $H_i$ under isomorphism map, are both equal to the binomial coefficient $\binom{p-1}{d}$.

(3) Check each time whether $Deck(G^d) = Deck(G)$ holds. If yes, then $G^d$ will be a proper extension of $G_i$ to $G$ i.e. $G^d \cong G$ (and hence $H^d$ will be a proper extension of $H_i$ to $H$ i.e. $H^d \cong H$). If not, then the extensions are improper and the added edges should be removed for trying the other choice.

(4) Since the deck of $G$ (and deck of $H$) is legitimate (**legitimate deck:** a deck which can actually be obtained from some graph) there must exist a choice for which $Deck(G^d) = Deck(G)$ implying $G^d \cong G$ and will hold. And $Deck(H^d) = Deck(H)$ Which will imply that $H^d \cong H$.

(5) We have achieved the isomorphic extension of $G_i$ to $G^d$ and $H_i$ to $H^d$, i.e. $G^d \cong H^d$, therefore, $G \cong H$.

□

**Proof 3:** Given unlabelled legitimate deck of point deleted subgraphs of $G$ we proceed to show that we can reconstruct a unique graph up to isomorphism (i.e. the reconstructed graph is isomorphic to $G$)

(1) Consider a subgraph in the deck.

(2) Add some d edges by joining some d points of this subgraph to a new point.

(3) Check subdeck property for this new graph and G.

(4) If it holds then stop. Else continue the steps from step (1) for all possible reconstructions equal to the binomial coefficient, $\binom{p-1}{d}$.

(5) At least one choice must exist for which the subdeck property will hold good and for the graph thus produced, say $G^*$, we will have $Deck(G^*) = Deck(G)$. Otherwise, the deck of $G$ will not be legitimate, a contradiction. If there exists only one choice (for adding d edges), then a unique graph is thus constructed and we are through. When there is more than one choice, we have to see that all these graphs are isomorphic. But this is clear since these multiple graphs are arrived at due to presence of more than one equivalent points in some equivalent set belonging to the equivalent partitioning of the point deleted subgraph in the deck taken for extension, whose **proper subset** is required to be used for extension.

□

**Remark 3.1:** If we carry out partitioning of $V(G_i)$ into equivalent sets $\{V_1, V_2, \cdots, V_k\}$ and let $\{|V_1|, |V_2|, \cdots, |V_k|\}$ be their cardinalities. Let $s_1, s_2, \cdots, s_k$ be the number of points in the sets $V_1, V_2, \cdots, V_k$ respectively which when joined to the external point produce a graph $G^+$ (a copy of G) such that $Deck(G^+) = Deck(G)$. It is clear that we can construct in all

$$\prod_{i=1}^{k} \binom{|V_i|}{s_i}$$

number of graphs isomorphic to $G^+$.

**4. A Counting Formula:** The count of isomorphic graphs that one will construct following the procedure discussed above in the proof 3 can be easily obtained as follows:

**Theorem 4.1:** Let $G_i$ be a graph in the deck of $G$. And let $\{V_1, V_2, \cdots, V_k\}$ be the partitioning of the vertex set of $G_i$ into equivalent sets, and $\{|V_1|, |V_2|, \cdots, |V_k|\}$ be their cardinalities. Let $\{r_1, r_2, \cdots, r_k\}$ be the points that are required to be joined by an edge to the new point to be taken outside to reconstruct $G$. Then the count of distinct isomorphic graphs that can be achieved through reconstruction is equal to the following product of the binomial coefficients:

$$\prod_{j=1}^{k} \binom{|V|_j}{r_j}$$

**Proof:** Suppose an equivalent set $V_j$ has to contribute $r_j$ points to join with the external point, then **any** $r_j$ among the $|V_j|$ points can be chosen (by definition). Also such choice for each equivalent set can be made independently. Hence the result. □

**5. The Reconstructibility of $S_Q^s(G)$:** Kelly lemma [3] shows the reconstructibility of $S_Q(G)$, the number of copies of subgraph Q in G, when $|V(G)| > |V(Q)|$. We now proceed to show that the number of copies of subgraph Q in G, $S_Q^s(G)$, is reconstructible, where $|V(G)|=|V(Q)|$.

**Theorem 5.1:** The number $S_Q^s(G)$ where Q is a spanning subgraph of G, i.e. when $|V(G)|=|V(Q)|$, is reconstructible.

**Proof:** Consider a graph $G_i \equiv G - v_i$ of G in the deck of G. Take a new point $v$ outside. Since $Deck(G)$ is legitimate G is reconstructible by theorem 3.2 and so there should exist points (at least equal to degree of $v_i$ in number) among (p – 1) points of $G_i$, such that any one point among which when we join to a new point $v$ taken outside and form an extended graph, say $G_{j_1}^1$, then it will be a subgraph of G. (when no such a graph can be formed then degree of $v_i$ will be zero, G will be disconnected, and so $S_Q^s(G) = 0$ hence we are done.). Let the number of copies of spanning subgraphs Q in $G_{j_1}^1$ be equal to $\alpha_{j_1}^1$ (*Note that each such a copy of Q contains the newly added edge.*). Since by Kelly lemma the degree list of G is reconstructible, so let $d(v) = d(v_i) = d$ say. So there will be **at least** d points $x_{j_1}, x_{j_2}, \cdots, x_{j_d} \in V(G_i)$ which when joined, one at a time, to the point taken outside, we will form the graphs $G_{j_1}^1, G_{j_2}^1, \ldots, G_{j_d}^1$ which will be subgraphs of G. Let number of copies of spanning subgraphs Q in these graphs be $\alpha_{j_1}^1, \alpha_{j_2}^1, \ldots, \alpha_{j_d}^1$.

Since $Deck(G)$ **is legitimate** there should exist pairs of points (when $2 \leq d$) among (p – 1) points of $G_i$ a pair among which when we join to a new point $v$ taken outside and form an extended graph $G_{j_1}^2$ then it will be a subgraph of G. Let the number of copies of spanning subgraphs Q in $G_{j_1}^2$ equal to $\alpha_{j_1}^2$ (*Note that here each such a copy of Q contains both the newly added edges.*). Actually there will be **at least** $\binom{d}{2}$ pairs ($2 \leq d$), made up of points $x_{j_1}, x_{j_2}, \cdots, x_{j_d} \in V(G_i)$, such that when these pairs are joined, one at a time, to the point taken outside we will form the graphs $G_{j_1}^2, G_{j_2}^2, \ldots, G_{j_{\binom{d}{2}}}^2$ which will be subgraphs of *G*. Let the number of copies of spanning subgraphs Q in these graphs be $\alpha^2_{j_1}, \alpha^2_{j_2}, ..., \alpha^2_{j_{\binom{d}{2}}}$.

Since $Deck(G)$ **is legitimate** there should exist triplets, **at least** $\binom{d}{3}$ in number, of points (when $3 \leq d$) among $(p-1)$ points of $G_i$ a triplet among which when we join to a new point $v$ taken outside and form an extended graph $G^3_{j_1}$ then it will be a subgraph of *G*. Continuing on similar lines let the number of copies of spanning subgraphs Q (*containing all the three newly added edges*) that we form in this case be $\alpha^3_{j_1}, \alpha^3_{j_2}, ..., \alpha^3_{j_{\binom{d}{3}}}$.

Continuing on similar lines there should exist **at least** one d-tuple made up of points $x_{j_1}, x_{j_2}, \cdots, x_{j_d} \in V(G_i)$ such that when all these points are joined to the point taken outside we will form the graph $G^d_{j_1}$, which will be isomorphic to *G* itself. Let the number of copies of spanning subgraphs Q (*containing all the d newly added edges*) that we form in this case be $\alpha^d_{j_1}$.

So, let there be **the** d points, $x_{j_1}, x_{j_2}, \cdots, x_{j_d} \in V(G_i)$ such that when all these points are joined to the external point we get the extended graph $G^d_{j_1}$ which is isomorphic to graph *G*, and let $\{\alpha^1_{j_1}, \alpha^1_{j_2}, ..., \alpha^1_{j_d}\}$, $\{\alpha^2_{j_1}, \alpha^2_{j_2}, ..., \alpha^2_{j_{\binom{d}{2}}}\}$, $\{\alpha^3_{j_1}, \alpha^3_{j_2}, ..., \alpha^3_{j_{\binom{d}{3}}}\}$, ... , $\{\alpha^d_{j_1}\}$ be the counts of graphs isomorphic to Q derived from 1-subsets, 2-subsets, 3-subsets, d-subset of $\{x_{j_1}, x_{j_2}, \cdots, x_{j_d}\}$ respectively, such that the **edges** added by joining the points of these r-subsets, r = 1, 2, **...**, d are **all present in the respective graphs isomorphic to Q** obtained from each of these r-subsets, then $S^s_Q(G)$ equal to sum of all these $\alpha^r_{j_s}$, s = 1, 2, ..., $\binom{d}{r}$.

Thus, $S^s_Q(G)$ is reconstructible. □

**6. Complexity of the Reconstruction Algorithm:** The complexity of the procedure described above of reconstructing the graph from its deck is dependent on the number of sets in the equivalent partitioning of the vertex set of a point deleted subgraph in the deck that we take for extension. Suppose the extension at a point of the subgraph in the deck chosen for extension leads to violation of subdeck property then the extension of that subgraph at any other point belonging to the same equivalent set in which the first point lies will also lead to violation of subdeck property. We can therefore discard that entire equivalent set for the extension purpose.

**Definition 6.1:** Let $G_i$ be a graph in the deck of $G$. And let $\{V_1, V_2, \cdots, V_k\}$ be the partitioning of the vertex set of $G_i$ into equivalent sets, and $\{|V_1|, |V_2|, \cdots, |V_k|\}$ be their cardinalities. A set $V_i$ is called **valid set** if for at least at one element (point) of this set the extension can be carried out maintaining the subdeck property.

**Definition 6.2:** A **set** $V_i$ is called **invalid set** when it is not valid.

**Theorem 6.1:** Given the partitioning of the vertex set of some graph $G_i$ in the deck of G, the order of complexity of reconstructing a graph from its deck is equal to $\sim O(T \times d(v_i) \times \Phi)$, where T is the time required for the checking of the subdeck property (described above) and $d(v_i)$ is degree of some vertex $v_i$ of G. And $\Phi$ is the number of invalid sets in the partitioning of the vertex set of $G_i$ into equivalent sets, for extension.

**Proof:** Take for extension purpose some subgraph $G_i \equiv G - v_i$ of G such that $d(v_i)$ is the degree $v_i$.
  (1) Take some equivalent set in the partitioning of the vertex set of $G_i$ into equivalent sets and carry out extension at some element (point) of it. Check the subdeck property, in time T say.
  (2) If the extension is valid, try to extend at some other elements of the same set, one by one, till one gets invalidity of subdeck property or till the entire set is consumed.
  (3) If the extension is invalid then discard that entire set from further consideration for extension and repeat step (1) starting with some other equivalent set.
  The proof is now clear. □

**Example:** Consider the following deck for reconstructing the graph [3].

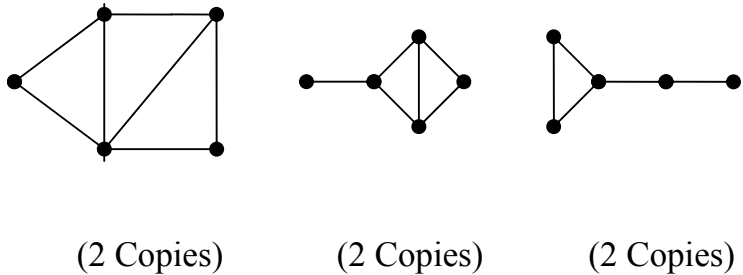

(2 Copies)    (2 Copies)    (2 Copies)

Proceeding as per the algorithm (discussed in **proofs of theorem 3.1**), we reconstruct the following graph from the given deck:

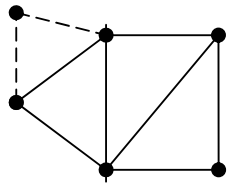

**7. The Edge-Reconstruction conjecture:** The most natural variation of the Kelly-Ulam conjecture is an analogue for deletion of edges. Let $G$ be a simple (n, e) graph and let the collection of all graphs $\{G - e_i\}$, $e_i \in E(G)$, the edge set of $G$ is given. The graph $G$ is said to be edge-reconstructible if it is uniquely determined by this deck. The edge-reconstruction conjecture due to Harary [4] states that: Every simple (n, e) graph, e > 3, is edge-reconstructible. Greenwell [5] has proved the following theorem:

**Theorem 7.1(Greenwell, 1971):** A reconstructible graph with no isolated vertices is edge-reconstructible.

It can be seen that:

**Theorem 7.2:** A graph with k, k > 0, isolated vertices is edge-reconstructible if and only if the graph obtained by deleting these k isolated vertices is edge-reconstructible.

**Proof:** Edge-Kelly lemma [3] permits the edge-reconstruction of the degree sequence and thus the isolated points. Hence every edge-reconstruction from a deck $\{G - e_i, i = 1, 2, \cdots, q$ and $e_i \in E(G)$, the edge set of $G$.$\}$ of a simple (p, q) graph $G$, has the same number of isolated points, which are all present in each graph in the deck.

□

Hence, in the light of the development made so far in this paper and the theorems 7.1, 7.2, we have the following

**Corollary 7.1 (Harary):** Every graph is edge-reconstructible.

**8. The Direct Approach for Edge Reconstruction:** By proceeding exactly on similar lines and developing the similar ideas for line (edge) deleted decks one can achieve the independent fixing of the edge-reconstruction conjecture as follows:

Edge-reconstruction conjecture (Harary [4]) can be stated as follows: Let $G$ and $H$ be two simple (p, q) graphs with p points and q lines (edges), $q > 3$, and let there exists a (1-1), onto map, called edge-hypomorphism:
$$\theta: E(G) \to E(H),$$
$$e_k \in E(G) \to f_k \in E(H)$$
such that $G_k \equiv G - e_k \cong H_k \equiv H - f_k, \forall k, k = 1, 2, \cdots, q$, then $G \cong H$. (The symbol $\cong$ stands for isomorphism.)

The existence of such map is called edge-hypomorphism and thus: the conjecture states that edge-hypomorphism implies isomorphism.

Let $e_i \in E(G)$ then the subgraph $G_i \equiv G - e_i$ is called the edge deleted subgraph (EDS) of $G$. The collection of all EDS, $\{G_i \equiv G - e_i, i = 1, 2, \cdots, q\}$, is called deck of (edge deleted subgraphs of) $G$, $Deck(G)$. The graph $G$ is said to be edge-reconstructible if every graph with the same deck as $G$ is isomorphic to $G$. Thus, the conjecture can be stated as: Every graph with more than three edges is uniquely (up to isomorphism) edge-reconstructible from its deck.

In order to reconstruct $G$ from its deck, $Deck(G)$, it suffices to identify which must be the vertex-pair of nonadjacent vertices for the missing edge $e_i$ in some graph $G_i$ in the deck. We proceed to accomplish this task.

Let $G$ and $H$ be two simple (p, q) graphs, i.e. having p points (vertices) and q lines (edges), and let $G$ be isomorphic to H. Then there will exist a (1-1), onto, adjacency preserving map
$$\theta: V(G) \to V(H),$$
$$v_i \in V(G) \to u_i \in V(H)$$
such that if pair of vertices, $v_k, v_l \in V(G)$, the vertex set for $G$, are adjacent (nonadjacent) then the corresponding pair of vertices, under the isomorphism map $\theta$, $u_k, u_l \in V(H)$, the vertex set for H, will be adjacent (nonadjacent).

Suppose now that the vertices, $v_k, v_l$ of G are nonadjacent and suppose we make them adjacent and also we add an edge in the edge set of H, E(H), by joining the points $u_k, u_l \in V(H)$ such that $u_k = \theta(v_k), u_l = \theta(v_l)$, and denote the extended G and H by $G^e$ and $H^e$ then clearly the supergraphs $G^e$ and $H^e$ will be isomorphic under the **same** isomorphism map.

Few definitions are now in order:

**Definition 8.1:** Let G and H be two isomorphic graphs. The process described above of extending graphs G and H to their respective supergraphs $G^e$ and $H^e$ is called the **isomorphic extension of isomorphic graphs.**

**Definition 8.2:** Let G be an unlabeled (p, q) graph and let $G^e$ be a supergraph of G obtained by adding a new edge, not in the edge set E(G), between some (unspecified) nonadjacent vertices of set V(G) is called the **extension** of G to $G^e$.

**Definition 8.3:** The subset of nonadjacent pairs of points (nonadjacent vertex-pairs) of E(G) is called a **set of equivalent vertex-pairs** if the extension of G by adding a new edge at any vertex-pair among these vertex-pairs, leads to graphs which are all isomorphic.

**Definition 8.4:** The collection of subsets $\{V_1, V_2, \cdots, V_m\}$ of **nonadjacent vertex-pairs** of the vertices in V(G), the vertex set of G, is called **partitioning of nonadjacent vertex-pairs of V(G) into equivalent sets** if all subsets $V_i$, i = 1, 2,…, m are equivalent sets,

$$V_i \cap V_j = \phi, \forall i \neq j, \text{ and } \bigcup_{i=1}^{m} V_i$$ represents the collection of all

the nonadjacent vertex-pairs, and where $\phi$ is a null set.

**Theorem 8.1:** Let G be a graph and let $\{V_1, V_2, \cdots, V_k\}$ be the partitioning of the of nonadjacent vertex-pairs of the vertices in V(G) of G into equivalent sets, and let $|V_1|, |V_2|, \cdots, |V_k|$ be their cardinalities. The number of nonisomorphic supergraphs that can be obtained from a graph $G$ by adding a single edge is equal to k, the total number of equivalent sets.

**Proof:** Any two graphs, one obtained by joining a nonadjacent vertex-pair belonging to some set, and the other obtained by joining some other nonadjacent vertex-pair belonging to the same set (though distinct ones) will be clearly isomorphic (refer to definition 2.3). Any two graphs, one obtained by joining nonadjacent vertex-pair belonging some set, and the other obtained by joining some other nonadjacent vertex-pair belonging to a different set will be clearly nonisomorphic.

□

**Definition 8.5:** Let $G$ and $H$ be two graphs on p points. We say that $Deck(G)$ **is a subdeck of** $Deck(H)$, or, $Deck(G) \subseteq Deck(H)$ when there exist a bijection $\psi : E(G) \to E(H)$ such that
$$e_i \in E(G) \to f_i \in E(H), \text{ where } E(G), E(H) \text{ are}$$
the edge sets of $G$ and $H$ respectively, and $G_i \equiv G - e_i$ is a subgraph of $H_i \equiv H - f_i$ for all i = 1, 2, …, p. Whenever $G_i \equiv G - e_i$ is a proper subgraph of $H_i \equiv H - f_i$ for at least one $i \in \{1, 2, …, p\}$, we say that $Deck(G)$ **is a proper subdeck of** $Deck(H)$, or, $Deck(G) \subset Deck(H)$.

Thus, $Deck(G)$ is **not a subdeck of** $Deck(H)$ if there exists no map $\psi$ (see definition 2.5), such that the inclusion relation, namely, $Deck(G) \subseteq Deck(H)$ holds good, i.e. this relation is violated for at least one subgraph of $G$ and the corresponding subgraph of $H$.

**Definition 8.6:** Checking whether $Deck(G) \subseteq Deck(H)$ holds (or does not hold) for some graphs $G$ and $H$ is called **checking of subdeck property for graphs $G$ and $H$**.

**Definition 8.7:** An **extension** of $G_i$ to $G^e$ is called **proper or valid** when the extension satisfies the subdeck property, viz.
$Deck(G^e) \subseteq Deck(G)$ when $G_i \in Deck(G).$

**Definition 8.8:** The subisomorphism map for two p point graphs $G$ and $H$ is a bijection
$$\psi : V(G) \to V(H) \text{ where}$$
$$v_i \in V(G) \to u_i \in V(H)$$

such that if a pair of vertices $(v_k, v_l) \in V(G)$ are adjacent then the corresponding pair of vertices, under the subisomorphism map $\theta$, $(u_k, u_l) \in V(H)$, will be adjacent and when $(v_k, v_l) \in V(G)$ are nonadjacent $(u_k, u_l) \in V(H)$ **may or may not be** adjacent.

**Remark 8.1:** Under existence of subisomorphism map $G$ will be a spanning subgraph of $H$.

**Lemma 8.1:** Let $G_1$ be a spanning subgraph of $G_2$, then $Deck(G_1) \subseteq Deck(G_2)$. Moreover, if $G_1$ is a proper subgraph of $G_2$ then $Deck(G_1) \subset Deck(G_2)$.

**Proof:** Since, $G_1$ is a subgraph of $G_2$ every adjacent vertex pair of $G_1$ will be adjacent in $G_2$ and $Deck(G_1) \subseteq Deck(G_2)$. When $G_1$ will be a proper subgraph of $G_2$ then the inclusion will be strict for at least one edge deleted subgraph of $G_1$ as a subgraph in the corresponding edge deleted subgraph of $G_2$. Since, there will exist a nonadjacent vertex pair in $G_1$ which will be adjacent in $G_2$ so when we consider an edge deleted subgraph of $G_1$ and the corresponding edge deleted subgraph of $G_2$ arrived at by deleting an edge that is present in both $G_1$ and $G_2$, $Deck(G_1) \subset Deck(G_2)$ becomes clear. In short for every edge in $G_1$ there is an edge in $G_2$ under the subisomorphism map. And for some nonadjacent pairs of vertices in $G_1$ some of the image pairs in $G_2$ may be adjacent. □

**Lemma 8.2:** Let $G_1$ and $G_2$ be two graphs on p points. If $Deck(G_1) \subseteq Deck(G_2)$ then $G_1$ is subgraph of $G_2$. Moreover, if $Deck(G_1) \subset Deck(G_2)$ then $G_1$ is a proper subgraph of $G_2$.

**Proof:** Suppose $G_1$ is not a subgraph of $G_2$ then there will not exist a subisomorphism map between their vertex sets. So, for every bijection between their vertex sets there will exist a pair of vertices $(v_k, v_l) \in V(G_1)$ which is adjacent while the corresponding pair of vertices $(u_k, u_l) \in V(G_2)$ (under the bijection) will be nonadjacent. So, if we form $Deck(G_1)$ and $Deck(G_2)$ and check under the same bijection whether the edge deleted subgraphs of $G_1$ are subgraphs of the corresponding edge deleted subgraphs of $G_2$ we will get violation for some subgraph obtained by deletion of some edge other than the one corresponding to edge joining $v_k, v_l$.

If we combine lemma 2.1 and lemma 2.2 we have

**Theorem 8.2:** Let $G_1$ and $G_2$ be two graphs on p points. Then $G_1$ be a subgraph of $G_2$ if and only if $Deck(G_1) \subseteq Deck(G_2)$. Moreover, $G_1$ is a proper subgraph of $G_2$ if and only if $Deck(G_1) \subset Deck(G_2)$.

**9. Edge-Reconstruction Conjecture:** We now proceed to show that edge-hypomorphic graphs are also isomorphic, i.e. the edge-reconstruction conjecture is true for all graphs.

**Theorem 9.1 (Edge-Reconstruction Conjecture):** Edge-Hypomorphic graphs are isomorphic.

**Proof 1:** Let $G$ and $H$ be two simple edge-hypomorphic graphs on q lines, q > 3, i.e. let there exists a bijection,
$$\theta : E(G) \to E(H),$$
$$e_i \in E(G) \to f_i \in E(H)$$
such that $G_i \equiv G - e_i \cong H_i \equiv H - f_i, \forall i, i = 1, 2, \cdots, q$.
Let $\omega_i : V(G_i) \to V(H_i)$, be the corresponding isomorphism maps for i = 1, 2, …, q.

(1) Consider some subgraph $G_i \equiv G - e_i$ of G and the corresponding isomorphic subgraph $H_i \equiv H - f_i$ of H.

(2) Extend isomorphically the graphs $G_i$ to $G^1$ and $H_i$ to $H^1$ by adding an edge by joining some nonadjacent vertex-pair $(x_j, x_k) \in E(G_i)$, the edge set of $G_i$, leading to the formation of the graph $G^1$ and by adding an edge joining the corresponding nonadjacent vertex-pair $(y_j, y_k) \in E(H_i), y_j = \omega_i(x_j), y_k = \omega_i(x_k)$, leading to the formation of the graph $H^1$. (Note that there are in all $\Phi$ number of possibilities for adding edges, one in $G_i$ and the corresponding one under isomorphism map in $H_i$, where $\Phi$ is the count of nonadjacent vertex-pairs in both $G_i$ and $H_i$.)

(3) Check whether $Deck(G^1) \subseteq Deck(G)$ holds. If yes, then $G^1$ will be a proper extension of $G_i$ towards G (and automatically $H^1$ will be a proper extension of $H_i$ towards $H$). If not, then the extensions are improper and the added edges should be removed.

(3.1) When the extension is proper, we are done, and we stop. It is clear that $Deck(G^1) = Deck(G)$ when the extension is proper, since $G_i$ is short of only one edge than $G$.

(3.2) When the extension is improper then remove the added edges, both in $G_i$ and $H_i$ and again go to step (2). Again, extend isomorphically the graphs $G_i$ to $G^1$ and $H_i$ to $H^1$ by adding an edge (other than the one previously added and found improper among the ($\Phi$-1) edges), leading to the formation of the graph $G^1$ and the corresponding graph $H^1$. Go to step (3) and check subdeck property for graphs $G^1$ and $G$.

(4) Continue extending isomorphically, and each time taking the proper action implied by the correctness (or incorrectness) of adding the edges as per the above given steps, till we reach the stop action in step (3.1). At this stage (completing step (3.1)) we will not only have $Deck(G^1) \subseteq Deck(G)$ and $Deck(H^1) \subseteq Deck(H)$ but we will achieve more. We will actually have $Deck(G^1) = Deck(G)$ and $Deck(H^1) = Deck(H)$. Thus, we have achieved the isomorphic extension of $G_i$ to $G^1$ and $H_i$ to $H^1$, i.e. $G^1 \cong H^1$.

(5) Note that if the isomorphic extension of $G_i$ fails at each addition, i.e. if we cannot add an edge in $G_i$ then clearly $Deck(G)$ will be illegitimate.

**Proof 2:** Given unlabelled legitimate deck of edge deleted subgraphs of $G$ we proceed to show that we can edge-reconstruct a unique graph up to isomorphism (i.e. the edge-reconstructed graph will be isomorphic to $G$)

(1) Take a subgraph in the deck.
(2) Add some an edge by adding an edge between some nonadjacent vertex-pair of this subgraph.
(3) Check subdeck property for this new graph and G.
(4) If it holds then stop. Else continue the steps from step (2) for all possible reconstructions, in all equal to $\Phi$, the total number of nonadjacent vertex-pairs in the subgraph under consideration.
(5) At least one choice must exist for which the subdeck property will hold good and for the graph thus produced, say $G^*$, we will have $Deck(G^*) = Deck(G)$. Otherwise, the deck of $G$ will not be legitimate, a contradiction. If there exists only one choice, then a unique graph is thus constructed and we are through. When there is more than one choice, we have to see that all these graphs are isomorphic. But this is clear since these multiple graphs are arrived at having the same deck due to presence of more than one equivalent vertex-pairs in equivalent set to which the vertex-pair, causing valid extension, belongs.

□

**10. A Simple Count:** The count of isomorphic graphs that one will construct from a given legitimate deck of edge deleted subgraphs following the procedure discussed above in the proof 2 can be easily obtained as follows:

**Theorem 10.1:** Let $G_i$ be a graph in the $Deck(G)$. The count of distinct isomorphic graphs that can be achieved through edge-reconstruction is equal to the cardinality of the equivalent sets of nonadjacent vertex pairs in the partitioning of nonadjacent vertex-pairs of the vertices in $V(G_i)$ into equivalent sets for extension.

**Proof:** Every equivalent set give rise to only one graph up to isomorphism by adding an edge.

□

**11. Complexity of the Edge-Reconstruction Algorithm:** The complexity of the procedure described above of edge-reconstructing the graph from its deck is dependent on the number of sets in the equivalent partitioning of the of nonadjacent vertex-pairs of an edge deleted subgraph in the deck that we take for extension. Suppose the extension by joining a nonadjacent vertex-pair of the subgraph in the deck chosen for extension leads to violation of subdeck property then the extension at any other nonadjacent vertex-pair of that subgraph belonging to the same equivalent set in which the first nonadjacent vertex-pair lies will also lead to violation of subdeck property. We can therefore discard that entire equivalent set for the extension purpose.

**Definition 11.1:** Let $G_i$ be a graph in the deck of $G$. And let $\{V_1, V_2, \cdots, V_k\}$ be the partitioning of the nonadjacent vertex-pairs of $G_i$ into equivalent sets, and $\{|V_1|, |V_2|, \cdots, |V_k|\}$ be their cardinalities. A set $V_i$ is called **valid set** if the extension, by making a nonadjacent vertex-pair adjacent, maintains the subdeck property.

**Definition 11.2:** A **set** $V_i$ is called **invalid set** when it is not valid.

**Theorem 11.1:** The order of complexity of edge-reconstructing a graph from its deck is equal to $\sim O(T \times \Phi)$, where T is the time required for the checking of the subdeck property and $\Phi$ is the number of sets in the partitioning of the vertex set of $G_i$ into equivalent sets, for extension.

**Proof:** Take for extension purpose some subgraph $G_i \equiv G - e_i$ of G.
 (1) Take some equivalent set and carry out extension at some element (a nonadjacent vertex-pair) of it. Check the subdeck property, in time T say.
 (2) If the extension is valid, stop.
 (3) If the extension is invalid then discard that entire set from further consideration for extension and repeat step (1) starting with some other equivalent set.
 The proof is now clear.
 □

**12. Complete Set of Invariants and Reconstruction:** One more approach to reconstruction problem is to devise a complete set of invariants that determines a graph uniquely up to isomorphism and use that complete set for reconstruction. We now proceed to discuss it in this section.
 An invariant of a graph G is a number associated with G which has same value for every graph isomorphic to G. For example, the number of points and the number of lines of a graph are clearly invariant numbers. A complete set of invariants is that set of numbers which completely determine a graph up to isomorphism. No decent complete set of invariants for a graph is known [6].
 We now propose to associate a vector having certain numbers as components with a graph and show that these components of the vector together form a complete set of invariants for that graph. Using this complete set one can see that Kelly-Ulam conjecture for graphs follows.

**13. Preliminaries:** Let G be a (p, q) labeled connected graph, i.e. a graph containing p points and q lines. We denote by $C_r(G)$ the number of labeled (r, r – 1) trees in G, r ≤ p, and let us denote by $C(G)$ the vector of tree counts in G, thus:

$$C(G) = (C_1(G), C_2(G), \ldots, C_r(G), \ldots, C_p(G)).$$

We thus have $C_1(G) = p$ and $C_2(G) = q$, the invariant numbers mentioned above.

**Theorem 13.1:** Let $G$ and $H$ be two labeled connected graphs. $G \cong H$ if and only if $C(G) = C(H)$.

**Proof:** When $G \cong H$ there exists an adjacency preserving bijection between the vertex sets of $G$ and $H$, so essentially we can keep $G$ on $H$ (or $H$ on $G$) such that one graph completely hides the other graph and so $C(G) = C(H)$ becomes clear.

Conversely, let $C(G) = C(H)$. We proceed to show that for every (k, k – 1) tree in $G$ there exists a (k, k – 1) tree in $H$ **isomorphic** to it for all k ≤ p which leads to the outcome that for every spanning tree in $G$ there exists a spanning tree in $H$ **isomorphic** to it (k = p case). Using the fact that a graph is **reconstructible** from all its spanning trees we settle the claim.

Since there is only one (1, 0), (2, 1) and (3, 2) tree up to isomorphism therefore for every (1, 0), (2, 1) and (3, 2) tree in $G$ there is a (1, 0), (2, 1) and (3, 2) tree in $H$ isomorphic to it. Note that there are two (4, 3) trees $T_1^4$ and $T_2^4$ up to isomorphism as shown below and let they be in quantity $\alpha_i^4$ and $\beta_i^4$, $i = 1, 2$ in $G$ and $H$ respectively:

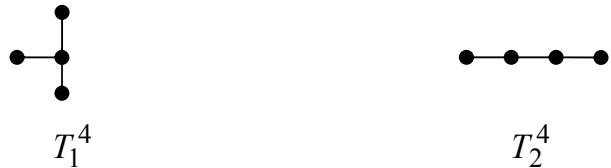

$T_1^4$ $\quad\quad\quad\quad\quad\quad T_2^4$

There are three types of (5, 4) trees $T_1^5$, $T_2^5$, $T_3^5$ up to isomorphism as shown below and let they be in quantity $\alpha_i^5$ and $\beta_i^5$, $i =$ 1 to 3 in $G$ and $H$ respectively:

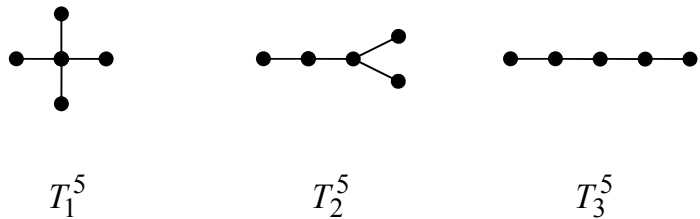

$T_1^5 \quad\quad\quad\quad T_2^5 \quad\quad\quad\quad T_3^5$

Let $T_1^k$ and $T_2^k$ be two (nonisomorphic) trees on $k$ points. Let $d(T_1^k)$ and $d(T_2^k)$ be their degree sequences written in nonincreasing order, namely,

$d(T_1^k) = d_1(T_1^k) \geq d_2(T_1^k) \geq d_3(T_1^k) \geq \ldots \geq d_k(T_1^k)$

$d(T_2^k) = d_1(T_2^k) \geq d_2(T_2^k) \geq d_3(T_2^k) \geq \ldots \geq d_k(T_2^k)$

**Definition 13.1:** $d(T_1^k) = d(T_2^k)$ if $d_j(T_1^k) = d_j(T_2^k)$ for all $j = 1$ to $k$ and $d(T_1^k) > d(T_2^k)$ if $d_j(T_1^k) = d_j(T_2^k)$ for $j = 1$ to $s$, $s < k$, then

$d_j(T_1^k) > d_j(T_2^k)$ for $j = (s+1)$.

**Definition 13.2:** $T_1^k > T_2^k$ if $d(T_1^k) > d(T_2^k)$ or when $d(T_1^k) = d(T_2^k)$ then the count of trees of type $T_u^{(k-1)}$ is more in $T_1^k$ than in $T_2^k$ while the count of trees of type $T_v^{(k-1)}$ is less in $T_1^k$ than $T_2^k$ when $T_u^{(k-1)} > T_v^{(k-1)}$.

Thus, the above definition is **recursive** and starts with $T_1^4 > T_2^4$.

Now, from equality $C(G) = C(H)$ we have

$$\alpha_1^4 + \alpha_2^4 = \beta_1^4 + \beta_2^4 \rightarrow (2.1)$$

so, we have

$$(\alpha_1^4 - \beta_1^4) + (\alpha_2^4 - \beta_2^4) = 0 \rightarrow (2.2)$$

Now, clearly $T_1^4 > T_2^4$. If $\alpha_1^4 > \beta_1^4$ then for equation (2.1) to hold we must have $\alpha_2^4 < \beta_2^4$. So, let $(\alpha_1^4 - \beta_1^4) = u$ then $(\alpha_2^4 - \beta_2^4) = -u$. Since a $T_1^4$ tree contains three (3, 2) trees while a $T_2^4$ tree contains two (3, 2) trees, therefore $G$ will contain u ($3u - 2u = u$) number of (3, 2) trees more than $H$, a contradiction. Thus, for every (4, 3) tree in $G$ there is a (4, 3) tree in $H$ isomorphic to it.

Similarly, we clearly have $T_1^5 > T_2^5 > T_3^5$. If $\alpha_1^5 > \beta_1^5$ then either $\alpha_2^5 < \beta_2^5$ or $\alpha_3^5 < \beta_3^5$ or both. It is easy to check from this that $G$ will contain **more** trees of type $T_1^4$ and **less** trees $T_2^4$ than $H$, a contradiction. We have thus shown that the equality of $C_j(G)$ and $C_j(H)$, $j = 1$ to $5$ implies that for every $(j, j-1)$ tree in $G$ there is a $(j, j-1)$ tree in $H$ isomorphic to it.

Now we proceed by induction and assume that the equality of $C_j(G)$ and $C_j(H), j = 1$ to $k$, implies that for every $(j, j-1)$ tree in $G$ there is a $(j, j-1)$ tree in $H$ **isomorphic** to it and proceed to show the same for $j = (k+1)$.

Let $T_j^{(k+1)}$, $j = 1$ to $l$, be the types of $(k+1, k)$ trees up to isomorphism and let $\alpha_j^{(k+1)}$ and $\beta_j^{(k+1)}$ be their count in the labeled graphs $G$ and $H$ respectively. Also, let $T_j^k$, $j = 1$ to $m$, be the types of $(k, k-1)$ trees up to isomorphism and suppose by **induction** that their **count** in the labeled graphs $G$ and $H$ is **equal** i.e. $\alpha_j^k = \beta_j^k$ for all $j = 1$ to $m$.

Also, we assume without loss of generality that
$$T_1^{(k+1)} > T_2^{(k+1)} > T_3^{(k+1)} > \ldots > T_l^{(k+1)}$$
and
$$T_1^k > T_2^k > T_3^k > \ldots > T_m^k$$
and our **aim** is to show that $\alpha_j^{(k+1)} = \beta_j^{(k+1)}$ for all $j = 1$ to $l$. For this we proceed by induction on the number of places $r$ where $\alpha_j^{(k+1)}$ and $\beta_j^{(k+1)}$ differ in their count. Now, if they differ in only one place i.e. if $r = 1$ then without loss of generality let $\alpha_j^{(k+1)} = \beta_j^{(k+1)}$ for all $j$ except $j = s$. But this is impossible because this will violate the equality of $C_{(k+1)}(G)$ and $C_{(k+1)}(H)$. So, let $r = 2$, i.e. let $\alpha_s^{(k+1)} > \beta_s^{(k+1)}$ and $\alpha_t^{(k+1)} < \beta_t^{(k+1)}$, for some $s < t$. But this will imply that there will exist trees $T_u^k$ and $T_v^k$, $T_u^k > T_v^k$, $u < v$, such that the count of trees $T_u^k$ as subtrees in $T_s^{(k+1)}$ will be more in $G$ than $H$, while the count of trees $T_v^k$ as subtrees in $T_t^{(k+1)}$ will be less in $G$ than $H$. This in turn implies that $\alpha_u^k > \beta_u^k$ while $\alpha_v^k < \beta_v^k$, a contradiction. Assuming that any $r = w$ terms different is impossible we can see that an additional differing term, to make $r = w + 1$, at upper or lower end can't avoid existence of two or more terms like $\alpha_u^k > \beta_u^k$ and $\alpha_v^k < \beta_v^k$ which is a contradiction. Hence, $\alpha_j^{(k+1)} = \beta_j^{(k+1)}$ for all $j = 1$ to $l$.

Thus, $C(G) = C(H)$ implies that for every labeled spanning tree of $G$ there will exist a labeled spanning tree in $H$

isomorphic to it. Now, it is a well-known fact that a graph is uniquely reconstructible from all its spanning trees, therefore
$G \cong H$.

□

We thus have the following

**Theorem 13.2:** The components of $C(G)$ forms a complete set of invariants for $G$.

**14. Kelly-Ulam Conjecture for (Connected) Graphs:** Let $G$ and $H$ be two simple connected (since the case for disconnected graphs is already settled [3]) graphs with p points, p>2, and let there exist a (1-1), onto map,
$\theta: V(G) \to V(H)$,
$v_i \in V(G) \to u_i \in V(H)$
such that $G_i \equiv G - v_i \cong H_i \equiv H - u_i, \forall i, i = 1,2,\cdots,p$, then
$G \cong H$. (The symbol $\cong$ stands for isomorphism.)
The existence of such a map is called hypomorphism and thus: the conjecture states that hypomorphism implies isomorphism.
   Let $v_i \in V(G)$. Then the induced subgraph $G_i \equiv G - v_i$ is called the vertex deleted subgraph (VDS) of $G$. The collection of all VDS, $\{G_i \equiv G - v_i, i = 1,2,\cdots,p\}$, is called the deck of (vertex deleted subgraphs of) $G$, $Deck(G)$. The graph $G$ is reconstructible if every graph with the same deck as $G$ is isomorphic to $G$. Thus, the conjecture can be stated as: Every graph with more than two vertices is uniquely (up to isomorphism) reconstructible from its deck.
   We proceed to show that in the light of theorems 13.1, 13. 2 it suffices to have the **knowledge of the degree of the missing vertex** $v_i$ for establishing the reconstructibility of $G$ **which we have** since the degree sequence of a graph is reconstructible [3].

**Theorem 15.1(Kelly-Ulam Conjecture):** Hypomorphic graphs are isomorphic.

**Proof:** Let $G$ and $H$ be two simple connected hypomorphic graphs on p points, p>2, i.e. let there exist a bijection,
$\theta: V(G) \to V(H)$,
$v_i \in V(G) \to u_i \in V(H)$, where
$G_i \equiv G - v_i \cong H_i \equiv H - u_i, \forall i, i = 1,2,\cdots,p$.

Let $\omega_i : V(G_i) \to V(H_i)$, be the corresponding isomorphism maps for i = 1, 2, ..., p.

Consider some subgraph $G_i \equiv G - v_i$ of G and the corresponding isomorphic subgraph $H_i \equiv H - u_i$ of H. If we will extend isomorphically the graphs $G_i$ to $G^1$ and $H_i$ to $H^1$ by adding an edge by joining some point $x_{j_1} \in V(G_i)$, the vertex set of $G_i$, to a new point $v$ not belonging to $V(G_i)$, leading to the formation of the graph $G^1$ and by adding an edge joining the corresponding point $y_{j_1} = \omega_i(x_{j_1}) \in V(H_i)$, the vertex set of $H_i$, to a new point $u$, not belonging to $V(H_i)$, leading to the formation of the graph $H^1$, then we have essentially formed two isomorphic spanning graphs.

Since $G_i \cong H_i$ we have for every subgraph of $G_i$ there is the corresponding isomorphic subgraph of $H_i$. Let $d(v_i) = d(u_i) = d$, the degrees of vertices $v_i$ and $u_i$ respectively. Now, we know that there are **some d neighbors** to $v_i$ and corresponding **d image points under isomorphism as neighbors** to $u_i$ in $G_i$ and $H_i$ respectively (which if we can specify then we achieve the reconstruction of G and H). **But whatever may be these d points**, once we **fix their choice**, among $\binom{(p-1)}{d}$ possibilities, then it is easy to see that:

(1) Every pair of spanning trees, one in $G_i$ and the other corresponding isomorphic one in $H_i$, will give rise to exactly d trees by isomorphic extension of these isomorphic trees counted in respectively $C_{(p-1)}(G_i)$ and $C_{(p-1)}(H_i)$. Therefore, we have $C_p^1(G) = d \times (C_{(p-1)}(G_i))$ and $C_p^1(H) = d \times (C_{(p-1)}(H_i))$, where $C_p^1(G)$ and $C_p^1(H)$ denotes **the count** of the spanning trees in G and H respectively such that the degree of $v_i$ and $u_i$ equal to **one.** It is **important to note at this juncture** that the spanning trees thus obtained may not be the spanning trees of G and H, but are definitely equal in count.

(2) Let us denote by $C_p^2(G)$ and $C_p^2(H)$ the **count** of the spanning trees that exist in respectively G and H such that the degrees of vertices $v_i$ and $u_i$ equal to **two**. We construct trees that are equal in count to trees in G and H as follows: we carry out all possible **disjoint** partitioning of the sets of (p-1) points of $G_i$ and $H_i$ into two sets $V^1(G_i)$ and $V^2(G_i)$ such that $V(G_i) = V^1(G_i) \cup V^2(G_i)$ and $V(H_i) = V^1(H_i) \cup V^2(H_i)$ such that if $x_{j_1} \in V^1(G_i)$ then $y_{j_1} = \omega_i(x_{j_1}) \in V^1(H_i)$ and if $x_{j_2} \in V^2(G_i)$ then $y_{j_2} = \omega_i(x_{j_2}) \in V^2(H_i)$. Let $d_1$ and $d_2$ be the points among the **d points chosen** in $G_i$ and the **images of these d points** under isomorphism that one gets in $H_i$, such that the $d_1$ points belong to $V^1(G_i)$ and their images to $V^1(H_i)$ while $d_2$ points belong to $V^2(G_i)$ and their images to $V^2(H_i)$ and $d_1$ and $d_2$ add up to d. It is easy to see that every subtree that spans the vertex set $V^1(G_i)$ and its isomorphic image that spans the vertex set $V^1(H_i)$ can be extended in $d_1$ ways by joining its point among the $d_1$ points to vertex $v_i$ and by joining the image point among the $d_1$ image points to vertex $u_i$. Similarly, every subtree that spans the vertex set $V^2(G_i)$ and its isomorphic image that spans the vertex set $V^2(H_i)$ can be extended in $d_2$ ways by joining its point among the $d_2$ points to vertex $v_i$ and by joining the image point among the $d_2$ image points to vertex $u_i$. Also, note that the choice of a point for extension among the points $d_1$ and $d_2$ can be made independently, so there will be in all **$d_1 \times d_2$ spanning trees** that can be obtained from **each pair of isomorphic trees**, one that spans the vertex set $V^1(G_i)$ and its isomorphic image that spans the vertex set $V^1(H_i)$ while the other that spans the vertex set $V^2(G_i)$ and its isomorphic image that spans the vertex set $V^2(H_i)$, by joining its point among the $d_1$ points to vertex $v_i$ and by joining the corresponding image point among the $d_1$ image points to vertex $u_i$ and by joining its point among the $d_2$ points to vertex $v_i$ and by joining the image point among the $d_2$ image points to vertex $u_i$.

Thus, by **considering**:
(a) All the possible disjoint partitioning of vertex set $V(G_i)$ and the corresponding automatic partitioning that follows by considering the images under isomorphism map in $V(H_i)$
(b) All the possible values of $d_1$ and $d_2$ offered by these partitioning of sets
(c) All the possible subtrees that span the partitioned sets
(d) All the possible extensions offered by each subtree and its isomorphic image

(e) All the possible spanning trees with degree of $v_i$ and $u_i$ equal to two and collecting them together we can have the **count** $C_p^2(G) = C_p^2(H)$.

Continuing with **disjoint** partitioning of $V(G_i)$ and the corresponding automatic partitioning that follows from isomorphism of $V(H_i)$ into **three subsets, four subsets, ..., d subsets**, etc. and forming all possible spanning trees with three, four, …, d, extension one from each subtree that spans the partitioned subset and from the point that is among the chosen d points etc. we can have the count of $C_p^3(G)$, $C_p^4(G)$, …, $C_p^d(G)$ equal to $C_p^3(H)$, $C_p^4(H)$, …, $C_p^d(H)$ respectively. Adding the respective counts, we thus have, $C_p(G) = C_p(H)$. Note that as far as the **count** is concerned the same outcome will result for any choice of d points.

Now, by **Kelly lemma** respectively $S_Q(G)$, the count of the copies of subgraph $Q$ in G is **reconstructible** when we are given a legitimate deck of graph G, *Deck*(G), and the cardinality of vertex sets $V(Q)$ is **less** than the cardinality of vertex set $V(G)$ [3]. Thus, from hypomorphism between G and H we are assured of the existence of an isomorphic subtree in H corresponding to every subtree of G when the cardinality of vertex set of that tree is **less** than the cardinality of vertex set $V(G)$. Thus, we have $C_j(G) = C_j(H)$, $j = 1$ to $(p-1)$. Hence $C(G) = C(H)$.

Thus, we have the **unique** components $C_j(G)$, $j = 1$ to $(p-1)$ by reconstructibility of $S_Q(G)$ and these unique components leads to a unique $C_p(G)$ so using theorem 2. 2 at this stage we have a unique G (up to isomorphism) from the given deck.

☐

**16. Another Complete Set of Invariants:** We now see that if we define a vector as follows by modifying the above considered invariant vector by replacing the components representing the count of the labeled trees of a particular size by count of the labeled trees with **multiplicity** i.e. a particular labeled (r, r − 1) tree is counted $k$ times if it occurs as a subtree of $k$ number of distinct (r+1, r) labeled trees that exist in G, then the components of this new vector also forms a complete set of invariants for G and by proceeding on similar lines we get the same development as follows:

Let G be a (p, q) labeled connected graph, i.e. a graph containing p points and q lines. Let us denote by $C*(G)$ the vector of tree counts in G, thus:

$$C*(G) = (C_1(G), C_2(G), ..., C_r(G), ..., C_p(G)).$$

where $C_p(G)$ denotes the count of all the labeled spanning trees in G, $C_{(p-1)}(G)$ denotes the count of (p-1, p-2) labeled trees. If these trees occur as subtrees more than once when we consider subtrees of every spanning tree separately, (i.e. as subtree of two or more distinct labeled spanning trees) then such trees are counted with that **multiplicity**. $C_{(p-2)}(G)$ is obtained from trees counted in $C_{(p-1)}(G)$ as subtrees with the **associated multiplicity**. Thus, trees counted in $C_r(G)$ are obtained as subtrees of $C_{(r+1)}(G)$ along with the associated multiplicity, for r = 1, 2, ..., p.

**Theorem 16.1:** Let G and H be two labeled connected graphs. $G \cong H$ if and only if $C*(G) = C*(H)$.

**Proof:** When $G \cong H$ there exists an adjacency preserving bijection between the vertex sets of G and H, so essentially we can keep G on H (or H on G) such that one graph completely hides the other graph and so $C(G) = C(H)$ becomes clear.

Conversely, let $C*(G) = C*(H)$. We proceed to show that for every spanning tree in G there exists a spanning tree in H **isomorphic** to it and using the fact that a graph is reconstructible from all its spanning trees we settle the claim.

Since there is only one (1, 0), (2, 1) and (3, 2) tree up to isomorphism therefore for every (1, 0), (2, 1) and (3, 2) tree in G there is a (1, 0), (2, 1) and (3, 2) tree in H isomorphic to it. There are two (4, 3) trees up to isomorphism, $T_1^4$ and $T_2^4$:

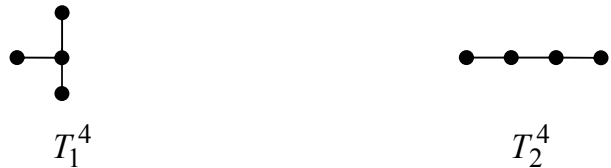

$T_1^4$  $T_2^4$

Using the equality of the respective components $C_j(G)$ and $C_j(H)$, j = 1 to 4, we show that for every labeled tree isomorphic to $T_1^4$ and $T_2^4$ in G there is a labeled tree isomorphic to $T_1^4$ and $T_2^4$ in H as follows: Let $\alpha_i^4$ and $\beta_i^4$, i = 1, 2 be the trees of type $T_1^4$ and $T_2^4$ in G and H respectively: we have

$$\alpha_1^4 + \alpha_2^4 = \beta_1^4 + \beta_2^4 \qquad \rightarrow (2.1)$$

so, we have

$$(\alpha_1^4 - \beta_1^4) + (\alpha_2^4 - \beta_2^4) = 0 \qquad \rightarrow (2.2)$$

Now it is easy to see that every (3, 2) tree is a subtree of a (4, 3) tree and a (4, 3) tree of type $T_1^4$ and of type $T_2^4$ contains respectively three and two (3, 2) trees as subtrees. Therefore, we have

$$3(\alpha_1^4 - \beta_1^4) + 2(\alpha_2^4 - \beta_2^4) = 0 \qquad \rightarrow (2.3)$$

The matrix of coefficients is

$$\begin{bmatrix} 1 & 1 \\ 3 & 2 \end{bmatrix}$$

which is of full rank, therefore, solving equations (2.2) and (2.3) simultaneously we get

$\alpha_1^4 = \beta_1^4$ and $\alpha_2^4 = \beta_2^4$. Thus, for every (4, 3) tree in $G$ there is a (4, 3) tree in $H$ isomorphic to it.

There are three types of (5, 4) trees $T_1^5$, $T_2^5$, $T_3^5$ up to isomorphism as shown below and let they be in quantity $\alpha_i^5$ and $\beta_i^5$, $i = 1$ to 3 in $G$ and $H$ respectively:

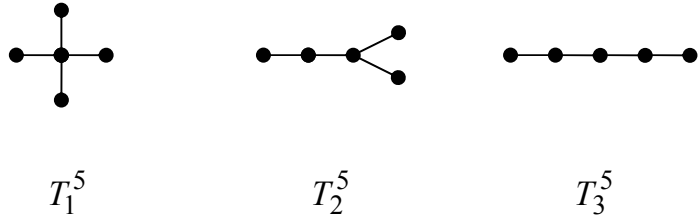

$$T_1^5 \qquad\qquad T_2^5 \qquad\qquad T_3^5$$

From equality of $C_j(G)$ and $C_j(H)$, $j = 1$ to 5 we have the following set of equations:

$$(\alpha_1^5 - \beta_1^5) + (\alpha_2^5 - \beta_2^5) + (\alpha_3^5 - \beta_3^5) = 0 \qquad \rightarrow (2.4)$$

Counting of trees of type $T_1^4$ as subtrees in these three trees we have

$$4(\alpha_1^5 - \beta_1^5) + (\alpha_2^5 - \beta_2^5) + 0 = 0 \qquad \rightarrow (2.5)$$

Counting of trees of type $T_2^4$ as subtrees in these three trees we have

$$0 + 2(\alpha_2^5 - \beta_2^5) + 2(\alpha_3^5 - \beta_3^5) = 0 \qquad \rightarrow (2.6)$$

The matrix of coefficients is

$$\begin{bmatrix} 1 & 1 & 1 \\ 4 & 1 & 0 \end{bmatrix}$$

which is of full rank, therefore, solving equations (2.4), (2.5) and (2.6) simultaneously we get
$\alpha_1^5 = \beta_1^5$, $\alpha_2^5 = \beta_2^5$ and $\alpha_3^5 = \beta_3^5$. Thus, for every (5, 4) tree in $G$ there is a (5, 4) tree in $H$ isomorphic to it.

Now proceed by induction and assume that the equality of $C_j(G)$ and $C_j(H)$, $j = 1$ to $k$, implies that for every $(j, j-1)$ tree in $G$ there is a $(j, j-1)$ tree in $H$ isomorphic to it and proceed to show the same for $j = (k+1)$. Let $T_j^{(k+1)}$, $j = 1$ to $l$, be the types of $(k+1, k)$ trees up to isomorphism and let $\alpha_j^{(k+1)}$ and $\beta_j^{(k+1)}$ be their count in the labeled graphs $G$ and $H$ respectively. Also, let $T_j^k$, $j = 1$ to $m$, be the types of $(k, k-1)$ trees up to isomorphism and suppose by **induction** that their count, counted with their multiplicities, in the labeled graphs $G$ and $H$ is equal i.e.
$\alpha_j^k = \beta_j^k$ for all $j = 1$ to $m$.

Also, we assume without loss of generality that
$$T_1^{(k+1)} > T_2^{(k+1)} > T_3^{(k+1)} > \ldots > T_l^{(k+1)}$$
and
$$T_1^k > T_2^k > T_3^k > \ldots > T_m^k$$
and our **aim** is to show that $\alpha_j^{(k+1)} = \beta_j^{(k+1)}$ for all $j = 1$ to $l$.
As is done above the first equation that we get (like equation (2.1) or (2.4)) will be
$$\sum_{j=1}^{l} (\alpha_j^{(k+1)} - \beta_j^{(k+1)}) = 0 \qquad \rightarrow (2.7)$$
which will lead to formation of the first row of the coefficient matrix containing all entries equal to **one**, as above.

Counting of trees of type $T_1^k$ as subtrees in these $l$ trees we will have the second row of the coefficient matrix with entries as nonnegative integers which form a **nonincreasing** sequence from left to right.

Counting of trees of type $T_2^k$ as subtrees in these $l$ trees we will have the third row of the coefficient matrix with entries as nonnegative integers which form a **nonincreasing** sequence **up to a certain term** and a **nonincreasing** sequence **afterwards** from left to right.

Counting of trees of type $T_3^k$ as subtrees in these $l$ trees we will have the third row of the coefficient matrix with entries as nonnegative integers which form a nonincreasing sequence up to a certain term (**greater than or equal to the term up to which the terms in the earlier row were nondecreasing**) and a nonincreasing sequence afterwards from left to right.

Continuing in this way up to last tree $T_m^k$, and counting of trees of $T_m^k$ type as subtrees in these $l$ trees we will have the last row of the coefficient matrix with entries as nonnegative integers which form a **nondecreasing** sequence from left to right.

Now, it can be seen that this **special** type of matrix with rows satisfying the above mentioned property is of **full rank**, which will imply that $\alpha_j^{(k+1)} = \beta_j^{(k+1)}$ for all $j = 1$ to $l$. This in turn implies that when $C^*(G) = C^*(H)$ then for every labeled spanning tree of $G$ there will exist a labeled spanning tree in $H$ isomorphic to it.

Now, it is a well-known fact that a graph is uniquely reconstructible from all its spanning trees, therefore
$$G \cong H. \quad \square$$

We thus have the following

**Theorem 16.2:** The components of $C^*(G)$ forms a complete set of invariants for $G$.

We proceed to show that in the light of theorem 16.2 it suffices to have the **knowledge of the degree of the missing vertex** $v_i$ for establishing the reconstructibility of $G$ **which we have** since the degree sequence of a graph is reconstructible [3].

**Theorem 16.3 (Kelly-Ulam Conjecture):** Hypomorphic graphs are isomorphic.

**Proof:** Let $G$ and $H$ be two simple connected hypomorphic graphs on $p$ points, $p > 2$, i.e. let there exist a bijection,
$$\theta : V(G) \to V(H),$$
$$v_i \in V(G) \to u_i \in V(H), \text{ where}$$
$$G_i \equiv G - v_i \cong H_i \equiv H - u_i, \forall i, i = 1, 2, \cdots, p.$$

Let $\omega_i : V(G_i) \to V(H_i)$, be the corresponding isomorphism maps for $i = 1, 2, \ldots, p$.

Consider some subgraph $G_i \equiv G - v_i$ of G and the corresponding isomorphic subgraph $H_i \equiv H - u_i$ of H. If we will extend isomorphically the graphs $G_i$ to $G^1$ and $H_i$ to $H^1$ by adding an edge by joining some point $x_{j_1} \in V(G_i)$, the vertex set of $G_i$, to a new point $v$ not belonging to $V(G_i)$, leading to the formation of the graph $G^1$ and by adding an edge joining the corresponding point $y_{j_1} = \omega_i(x_{j_1}) \in V(H_i)$, the vertex set of $H_i$, to a new point $u$, not belonging to $V(H_i)$, leading to the formation of the graph $H^1$, then we have essentially formed two isomorohic spanning graphs.

Since $G_i \cong H_i$ we have for every subgraph of $G_i$ there is the corresponding isomorphic subgraph of $H_i$. Let $d(v_i) = d(u_i)$ = d, the degrees of vertices $v_i$ and $u_i$ respectively. Now, we know that there are **some d neighbors** to $v_i$ and $u_i$ in $G_i$ and $H_i$ respectively (which if we can specify then we achieve the reconstruction of G and H). **But whatever may be these d points**, once we **fix** their choice, among $\binom{(p-1)}{d}$ possibilities, then it is easy to see that:

(1) Every pair of spanning trees, one in $G_i$ and the other corresponding isomorphic one in $H_i$, will give rise to exactly d trees by isomorphic extension of these isomorphic trees counted in respectively $C_{(p-1)}(G_i)$ and $C_{(p-1)}(H_i)$. Therefore, we have $C_p^1(G) = d \times (C_{(p-1)}(G_i))$ and $C_p^1(H) = d \times (C_{(p-1)}(H_i))$, where $C_p^1(G)$ and $C_p^1(H)$ denotes **the count** of the spanning trees **that can exist** in G and H respectively such that the degree of $v_i$ and $u_i$ equal to **one**.

(2) Let us denote by $C_p^2(G)$ and $C_p^2(H)$ the **count** of the spanning trees **that can exist** in respectively G and H such that the degrees of vertices $v_i$ and $u_i$ equal to two. We construct these trees as follows: we carry out all possible disjoint partitioning of the sets of (p-1) points of $G_i$ and $H_i$ into two sets $V^1(G_i)$ and $V^2(G_i)$ such that $V(G_i) = V^1(G_i) \cup V^2(G_i)$ and $V(H_i) = V^1(H_i) \cup V^2(H_i)$ such that if $x_{j_1} \in V^1(G_i)$ then $y_{j_1} = \omega_i(x_{j_1}) \in V^1(H_i)$ and if $x_{j_2} \in V^2(G_i)$ then $y_{j_2} = \omega_i(x_{j_2}) \in V^2(H_i)$. Let $d_1$ and $d_2$ be the points among the **d** **points chosen** in $G_i$ and the **images of these d points** under isomorphism one gets in $H_i$, such that the $d_1$ points belong to $V^1(G_i)$ and their images to $V^1(H_i)$ while $d_2$ points belong to $V^2(G_i)$ and their images to $V^2(H_i)$ and $d_1$ and $d_2$ add up to d. It is easy to see that every subtree that spans the vertex set $V^1(G_i)$ and its isomorphic image that spans the vertex set $V^1(H_i)$ can be extended in $d_1$ ways by joining its point among the $d_1$ points to vertex $v_i$ and by joining the image point among the $d_1$ image points to vertex $u_i$. Similarly, every subtree that spans the vertex set $V^2(G_i)$ and its isomorphic image that spans the vertex set $V^2(H_i)$ can be extended in $d_2$ ways by joining its point among the $d_2$ points to vertex $v_i$ and by joining the image point among the $d_2$ image points to vertex $u_i$. The choice of a point for extension among the points $d_1$ and $d_2$ can be made independently, so there will be in all $d_1 \times d_2$ spanning trees that can be obtained from each pair of trees, one that spans the vertex set $V^1(G_i)$ and its isomorphic image that spans the vertex set $V^1(H_i)$ while the other that spans the vertex set $V^2(G_i)$ and its isomorphic image that spans the vertex set $V^2(H_i)$, by joining its point among the $d_1$ points to vertex $v_i$ and by joining the image point among the $d_1$ image points to vertex $u_i$ and by joining its point among the $d_2$ points to vertex $v_i$ and by joining the image point among the $d_2$ image points to vertex $u_i$. Thus by considering:
(a) All the possible disjoint partitioning of vertex set $V(G_i)$ and $V(H_i)$
(b) All the possible values of $d_1$ and $d_2$ offered by these partitioning of sets
(c) All the possible subtrees that span the partitioned sets
(d) All the possible extensions offered by each subtree and its image
(e) All the possible spanning trees with degree of $v_i$ and $u_i$ equal to two
and collecting them together we can have the count $C_p^2(G) = C_p^2(H)$.

Continuing with disjoint partitioning of $V(G_i)$ and $V(H_i)$ into three subsets, four subsets, ..., d subsets, etc. and forming all possible spanning trees with three, four, ..., d, extension one from each subtree that spans the partitioned subset and from the point that is among the chosen d points etc. we can have the count of $C_p^3(G)$, $C_p^4(G)$, ..., $C_p^d(G)$ equal to $C_p^3(H)$, $C_p^4(H)$, ..., $C_p^d(H)$ respectively. Adding the respective counts, we thus have, $C_p(G) = C_p(H)$. Note that as far as the **count** is concerned this will happen for any choice of d points.

Now, by **Kelly lemma** respectively $S_Q(G)$, the count of the copies of subgraph Q is reconstructible when we are given a legitimate deck of graph G, *Deck(G)* such that the cardinality of vertex sets $V(Q)$ is **less** than the cardinality of vertex set $V(G)$ [3]. Thus, from hypomorphism between G and H we are assured of the existence of an isomorphic subtree in H corresponding to every subtree of G when the cardinality of vertex set of that tree is **less** than the cardinality of vertex set $V(G)$. Thus, counting all subtrees **with their multiplicities** we have $C_j(G) = C_j(H), j$ = 1 to $(p-1)$. Hence $C*(G) = C*(H)$.

Thus, we have the **unique** components $C_j(G), j = 1$ to $(p-1)$ by reconstructibility of $S_Q(G)$ and these unique components leads to a unique $C_p(G)$ so using theorem 2.2 at this stage we have a unique G (up to isomorphism) from the given deck.

$\square$

**17. Graphs reconstruction using reconstruction result for Trees:** In this section we begin with an elaboration of the Harary's formulation of reconstruction conjecture. This elaboration aims at exactly pointing out what is to be achieved for getting a complete solution to this problem. We then proceed with development of a new approach that uses the fact that every connected graph is uniquely reconstructible from its spanning trees, and every tree can be uniquely reconstructible from its pendant point deleted deck of trees [7].

**17.1 A light on Harary's Approach to Reconstruction Problem:** A graph property or parameter is reconstructible if its truth value or value (of that parameter) for every graph G is uniquely determined by its deck. It is well-known fact that from a given deck of one point deleted graphs of a graph the number of vertices, the number of edges, degree sequence, the number of subgraphs isomorphic to a given graph Q, when $|V(Q)| < |V(G)|$ etc. for the graph G are reconstructible. There are in fact some more other properties and parameters that are reconstructible.

Our elaboration of Harary's formulation of reconstruction conjecture goes as follows:

Let $\{G_i \cong G - v_i\}$ be a legitimate deck for some graph G.

(1) Choose and determine some one or more properties and/or parameters that are uniquely reconstructible from given deck.
(2) Carry out the constructions of all possible nonisomorphic graphs by extension, i.e. by taking a new vertex outside and joining this vertex by new edges connecting certain vertices of first graph, graph $G_1$ in the given deck of graphs such that all these extended nonisomorphic graphs satisfy the properties and parameters chosen in (1) and collect these graphs obtained by extension in set $S_1$.
(3) Carry out the constructions of graphs by extension satisfying the chosen properties and parameters as is done in (2) for every graph $\{G_i \cong G - v_i\}$ in the deck and collect them in respective sets $S_i$.
(4) Find set $S = \bigcap_i S_i$.

Reconstruction conjecture demands to show that set $S = \bigcap_i S_i$ is a **singleton** set.

So, using one or more properties and/or parameters that are uniquely reconstructible from given deck, if one can show that among the feasible graphs that can be obtained by extension from each member subgraph of the legitimate deck when we take their intersection only some unique graph survives than one has proved the reconstruction conjecture in the affirmative!

**17.2 Spanning Trees and Reconstruction:** In this sub section we give a new crisp proof of reconstruction conjecture! All graphs in this section are connected graphs.

We have based our proof upon two well known results:

*1*: A graph can be uniquely reconstructed from all its spanning trees.

*2*: A tree can be uniquely reconstructed from its pendant point deleted deck of subtrees.

**Theorem 17.2.1 (Reconstruction Conjecture):** Let $G$ and $H$ be two simple graphs with p points, p>2, and let there exist a (1-1), onto map,
$\theta: V(G) \to V(H)$,
$v_i \in V(G) \to u_i \in V(H)$
such that $G_i \equiv G - v_i \cong H_i \equiv H - u_i, \forall i, i = 1,2,\cdots, p$, then
$G \cong H$. (The symbol $\cong$ stands for isomorphism.)

The existence of such a map is called hypomorphism. The reconstruction conjecture states that hypomorphism implies isomorphism.

**Proof:** Let $\{T_1^G, T_2^G, \cdots\}$ and $\{T_1^H, T_2^H, \cdots\}$ be the spanning trees of $G$ and $H$ respectively.

Consider all spanning trees of $\{G_i \cong G-v_i\}$ and $\{H_i \cong H-u_i\}$.

Now, since $G_i \equiv G-v_i \cong H_i \equiv H-u_i, \forall i, i=1,2,\cdots,p$, therefore for every spanning tree of $\{G_i \cong G-v_i\}$ there will exists spanning tree of $\{H_i \cong H-u_i\}$ isomorphic to it. So, let $\{T_1^{v_i}, T_2^{v_i}, \cdots\}$ be the spanning trees of $\{G_i \cong G-v_i\}$ and $\{T_1^{u_i}, T_2^{u_i}, \cdots\}$ be the spanning trees of $\{H_i \cong H-u_i\}$ such that $T_j^{v_i} \cong T_j^{u_i} \forall j$.

Now, consider first tree $T_1^G$ in the above mentioned set of spanning trees of $G$, namely, $\{T_1^G, T_2^G, \cdots\}$. Let $\{v_{j_1}, v_{j_1}, \cdots v_{j_k}\}$ be the pendant points of $T_1^G$. In the spanning trees of $\{G_{j_l} \cong G-v_{j_l}, 1 \leq l \leq k\}$ there will exist suitable spanning trees which will give rise to tree $T_1^G$ by the well known result (*2*) mentioned above. Now, the isomorphic copies of these suitable spanning trees will exists in $\{H_{j_l} \cong H-u_{j_l}, 1 \leq l \leq k\}$ and they will give rise, again due to result (*2*), to a spanning tree, without any loss of generality say $T_1^H$, such that $T_1^G$ will be isomorphic to $T_1^H$. In this way we can build all spanning trees of $G$ and for every spanning tree of $G$ there will be exactly one isomorphic spanning tree of $H$.

Now, as per result (*1*) all spanning trees of $G$ will reconstruct unique $G$ and all spanning trees of $H$ will reconstruct unique $H$, therefore $G \cong H$ as desired!

□

### Acknowledgements

I am thankful to Prof. Dr. M. R. Modak, Bhaskaracharya Pratishthans, Pune, and Prof. Dr. T. T. Raghunathan, University of Pune, for their keen interest.